\documentclass[centertags,12pt]{amsart}
\usepackage{latexsym}
\usepackage{amssymb}
\numberwithin{equation}{section}
\makeatletter
\renewcommand{\subsection}{\@startsection
{subsection}{2}{0mm}{\baselineskip}{-0.25cm}
{\normalfont\normalsize\bf}}
\makeatother
\newtheorem*{proposition*}{Proposition}
\newtheorem*{corollary*}{Corollary}
\newtheorem{thm}{Theorem}

\newtheorem{fact}{Fact}
\newtheorem{theorem}{Theorem}[section]
\newtheorem{proposition}[theorem]{Proposition}
\newtheorem{lemma}[theorem]{Lemma}
\newtheorem{corollary}[theorem]{Corollary}
\theoremstyle{remark}

\newtheorem*{remark*}{Remark}
\newtheorem*{remarks*}{Remarks}
\newtheorem*{claim*}{Claim}
\def\F{\mathbb F}
\def\P{\mathbb P}
\def\N{\mathbb N}
\def\Z{\mathbb Z}
\def\cD{\mathcal D}
\def\cJ{\mathcal J}

\def\cO{\mathcal O}
\def\cK{\mathcal K}
\def\cW{\mathcal W}
\def\l{\ell}
\def\fq{\mathbb F_q}
\def\fl{\mathbb F_{\ell^2}}
\def\supp{{\rm supp}}
\def\df{{\rm div}}
\def\d0{{\rm div}_0}
\def\dinf{{\rm div}_{\infty}}
\def\dim{{\rm dim}}
\def\deg{{\rm deg}}
\def\sq{\sqrt{q}}
\def\frj{{\rm Fr}_{\mathcal J}}
\def\frx{{\rm Fr}_{X}}
\begin{document}
\author[R.~Fuhrmann]{Rainer Fuhrmann}
\author[F.~Torres]{Fernando Torres}\thanks{The paper was partially
written while Torres was visiting the
University of Essen (supported by the Graduiertenkolleg ``Theoretische
und experimentelle Methoden der Reinen Mathematik"); ICTP, Trieste
-Italy (supported by ICTP) and IMPA, Rio de Janeiro - Brazil (supported by
Cnpq)}
\title[Weierstrass points and optimal curves]{On Weierstrass points and 
optimal curves}
\address{Universit\"at
GH Essen, FB 6 Mathematik und Informatik, D-45117 Essen, Germany}
\address{IMECC-UNICAMP, Campinas-13083-970-SP, Brazil}
\email{RAINER.FUHRMANN@ZENTRALE.DEUTSCHE-BANK.dbp.de}
\email{ftorres@ime.unicamp.br}
\begin{abstract}
We use Weierstrass Point Theory and Frobenius orders  to prove the
uniqueness (up to isomorphism) of some optimal curves.
\end{abstract}
\maketitle
This paper continues the study, begun in \cite{ft} and \cite{fgt}, of 
curves over finite fields with many rational points, based on 
St\"ohr-Voloch's approach \cite{sv} to the Hasse-Weil bound by way of 
Weierstrass Point Theory and Frobenius orders. Some of the results were
announced in \cite{t}. 

A projective geometrically irreducible non-singular algebraic curve 
$X\! \mid\! \fq$ of genus $g$ is said to be {\it optimal} if
$$
\# X(\fq)= {\rm max}\{\# Y(\fq): \text{$Y\! \mid\! \fq$ curve of genus
$g$}\}\, .
$$ 
Optimal curves occupy a distinguished niche, for example, in coding theory 
after Goppa's \cite{go}. We recall that $\# X(\fq)$ is bounded from above
by the Hasse-Weil bound, namely
$$
q+2g\sq +1\, .
$$
The main goal of this paper is to sharpen and generalize some results in 
\cite{fgt}. In that paper Garcia and us improved and generalized previous
results obtained by R\"uck-Stichtenoth's \cite{rsti} and
Stichtenoth-Xing's \cite{stix}. 

We will mainly concerned with the uniqueness (up to
$\fq$-isomorphism) of some optimal curves. Roughly speaking, firstly the
Zeta Function of the curve is used to define a linear system on the curve.
Then, St\"ohr-Voloch's \cite{sv} is used to obtain the desired property.

We distinguish two cases according as $q$ is a square or not. 
In the first case, say $q=\l^2$, we look for curves $X\!\mid\!\fl$ that
attain the Hasse-Weil bound, that is, the so-called maximal curves. 
These curves were studied in \cite{stix}, 
\cite{gvl} (see also the references therein), \cite{ft} and \cite{fgt}. 
$X$ is equipped with the linear system $\cD:=|(\l+1)P_0|$,
$P_0\in X(\fl)$ (cf. \cite[\S1]{fgt}, \S2 here) and from the application
of 
\cite{sv} to $\cD$, Garcia and us deduced properties on the genus
and the uniqueness (up to $\fl$-isomorphism) of $X$ for some values of 
genus (cf. \cite{ft}, [op. cit.]). The following theorem follows 
from \cite{ft} and Proposition \ref{p2.1}; it improves \cite[Thm.\!
3.1]{fgt} and is a typical example of 
the results obtained here. We recall that the biggest genus that $X$ can 
have is $\l(\l-1)/2$ (cf. Ihara's \cite{ih}). 
\begin{thm}\label{A} Let $X\!\mid\! \fl$ be a maximal curve of genus $g$
and $\l$ odd.\newline
If $g>(\l-1)(\l-2)/4$, then
\begin{enumerate}
\item $X$ is $\fl$-isomorphic to the Hermitian curve $y^\l+y=x^{\l+1}$ so
that $g=(\l-1)\l/2$\, or
\item $X$ is $\fl$-isomorphic to the plane curve $y^\l+y=x^{(\l+1)/2}$ so 
that $g=(\l-1)^2/4$\, .
\end{enumerate}
\end{thm}
Theorem \ref{A}(1) is valid without restricting the parity of $\l$ if 
$g>(\l-1)^2/4$. Indeed, several characterizations of Hermitian curves
have already been given, see for example \cite{hstv} (and the references
therein), \cite{rsti} and \cite{ft}. See also Theorem \ref{t2.1}.

Furthermore, we show that the morphism associated to $\cD$ is an 
embedding. Hence we improve \cite[Prop.\!
1.10]{fgt}, that is, we can compute the genus of maximal 
curves under a hypothesis on non-gaps at $\fl$-rational points (see
\S\ref{2.3}).

Now we discuss the case $\sq\not\in \N$. Besides some curves of 
small genus, see for example \cite{carb}, the only known examples of
optimal curves, in this case, are the
Deligne-Lusztig curves associated to the Suzuki group $Sz(q)$ and to the
Ree group $R(q)$ \cite[\S11]{dl}, \cite{h}. They were studied in 
\cite{hsti}, \cite{h}, 
\cite{p} and \cite{hp}. Hansen and Pedersen \cite[Thm.\! 1]{hp} stated the
uniqueness, up to $\fq$-isomorphism, of the curve corresponding to $R(q)$
based on the genus, the number of $\fq$-rational points, and the group of
$\fq$-automorphisms of the curve.  They observed a similar result for the
curve corresponding to $Sz(q)$ (cf. \cite[p.100]{hp}) as a consequence of
its uniqueness up to $\bar \fq$-isomorphism (cf. \cite{he}). Hence, after
\cite{he} and 
\cite{hsti}, the curve under study in \S3 of this paper is 
$\fq$-isomorphic to 
the plane curve given by
$$
y^q-y=x^{q_0}(x^q-x)\, ,
$$
where $q_0=2^s$ and $q=2q_0$. This curve is equipped with the linear
system $g^4_{q+2q_0+1}=|(q+2q_0+1)P_0|$, $P_0$ a $\fq$-rational point.
By applying \cite{sv} to this linear system we prove (see \S3) the
\begin{thm}\label{B} Let $q_0, q$ be as above, $X\!\mid\!\fq$ a curve of 
genus $g$ such that:\quad (1)\, $g=q_0(q-1)$\quad and\quad
(2)\, $\#X(\fq)=q^2+1$. 

Then $X$ is $\fq$-isomorphic to the Deligne-Lusztig curve associated to 
the Suzuki group $Sz(q)$.
\end{thm}
We remark that a Hermitian curve can be also realized as a Deligne-Lusztig
curve associated to a projective special linear group (cf. \cite{h}). Then
its uniqueness (up to $\fq$) is also a consequence of its uniqueness up to
$\bar \fq$ (cf. \cite[p.100]{hp}).

A. Cossidente brought to our attention a relation between the curve in
Theorem \ref{B} and the Suzuki-Tits ovoid. This is described in the
Appendix.

It is our pleasure to thank: A. Garcia, R. Pellikaan and H. Stichtenoth
for useful conversations; A. Cossidente for having let us include his
observation in the Appendix. In addition, we want to thank Prof. 
J.F. Voloch for his interest in this work. 

{\bf Convention:} Throughout this paper by a curve we mean a projective
geometrically irreducible non-singular algebraic curve.
\section{Preliminaries}\label{1}
In this section we summarize some background material concerning
Weierstrass Point Theory, Frobenius orders and a rational divisor arising
from the Zeta Function of a curve defined over a finite field. 
\subsection{Weierstrass Point Theory}\label{1.1} 
Here we repeat relevant material from St\"ohr-Voloch's \cite[\S1]{sv} (see
also \cite[III.5]{fk}, \cite{lak} and \cite{sch}).

Let $X$ be a curve of genus $g$ over an algebraically closed field $k$ 
and $k(X)$ the field of 
rational functions on $X$. Let $\cD$ be a $g^r_d$ on $X$, say 
$$
\cD=\{E+\df(f): f\in \cD'\setminus\{0\}\}\subseteq |E|\, ,
$$
$E$ being an effective divisor on $X$ with $\deg(E)=d$ and $\cD'$ an 
$(r+1)$-dimensional $k$-subspace of $L(E):=\{f\in k(X)^*: 
E+\df(f)\succeq 0\}$. 

To any $P\in X$ one then associates the sequence of 
{\it $(\cD,P)$-orders}
$$
\{j_0(P)<\ldots<j_r(P)\}:=\{v_P(D): D\in \cD\}\, ,
$$
and on $X$ one has the so-called {\it $\cD$-ramification divisor}, namely 
$$
R=R^{\cD}=\df({\rm
det}((D^{\epsilon_i}_xf_j)))+\sum_{i=0}^{r}\epsilon_i\df(dx)+(r+1)E\, ,
$$
where $x$ is a separating variable of $k(X)|k$, $D^{i}_x$ is the
$i$th Hasse derivative, $f_0,\ldots,f_r$ is a $k$-base of $\cD'$, and
$(\epsilon_0,\ldots,\epsilon_r)$ (called {\it the $\cD$-orders}) is 
the minimum in the lexicographic order 
of the set
$$
\{(j_0,\ldots,j_r)\in \N^{r+1}: j_0<\ldots<j_r;\ {\rm
det}((D^{j_i}_xf_j))\neq 0\}\, .
$$
One has
\begin{equation}\label{eq1.1}
\begin{split}
& (a)\quad \deg(R)=\sum_{P\in
X}v_P(R)=\sum_{i=0}^{r}\epsilon_i(2g-2)+(r+1)d\, ,\\
& (b)\quad j_i(P)\ge \epsilon_i\qquad \text{for each $i$ and for each
$P$}\, ,\\
& (c)\quad v_P(R)\ge \sum_i(j_i(P)-\epsilon_i)\, ,\qquad \text{and} \\
& (d)\quad v_P(R)=\sum_i(j_i(P)-\epsilon_i)
\Leftrightarrow 
{\rm det}(\binom{j_i(P)}{\epsilon_j})\not\equiv 0 \pmod{{\rm
char}(k)}\, .\\
\end{split}
\end{equation}
Consequently the $(\cD,P)$-orders are  
$\epsilon_0,\ldots,\epsilon_r$ if and only if $P\in X\setminus \supp(R)$. 
The points in $\supp(R)$ are 
the so-called {\it $\cD$-Weierstrass points}.

The $\cK$-Weierstrass points, being $\cK=\cK_X$ the canonical linear 
system on $X$, are the {\it Weierstrass points} of $X$. In this case 
$H(P):= \N\setminus \{j_0(P)+1,\ldots,j_{g-1}(P)+1\}$ is the {\it 
Weierstrass semigroup} at $P$. We write  
$H(P)=\{m_0(P)=0<m_1(P)<\ldots\}$, the element $m_i(P)$ being the {\it
$i$th
non-gap} at $P$. The curve is called {\it classical} iff the $\cK$-orders
are $0,1,\ldots,g-1$ (i.e. $H(P)=\{0,g+1,g+2,\ldots\}$ for
$P\not\in\supp(R^{\cK})$).

To any $P\in X$ one also associates the {\it $i$th osculating plane} 
$L_i(P)$:  
via the identification $\cD\cong \P(\cD')^*$ each hyperplane $H$ in
$\P(\cD')$ correspond 
to a divisor $D_H\in \cD$; then $L_i(P)$ is the intersection of the
hyperplanes $H$ in $\P(\cD)^*$ such that $v_P(D_H)\ge j_{i+1}(P)$. Its 
(projective) dimension is $i$.  In terms of projective 
coordinates $L_i(P)$ can be described as follows:   
let $f_0,\ldots,f_r$ be a $(\cD,P)$-hermitian base of $\cD$, i.e. a 
$k$-base of $\cD'$ such that $v_P(t^{v_P(E)}f_i)=j_i(P)$ for
$i=0,\ldots,r$; $t$ being a local parameter at $P$. Then for
$i=0,\ldots,r-1$
\begin{equation}\label{eq1.2}
L_i(P)=H_{i+1}\cap\ldots\cap H_{r}\quad \mbox{with}\
D_{H_j}=E+\df(f_j),\ j=i+1,\ldots,r\, .
\end{equation}
\subsection{Frobenius orders}\label{1.2}
In the remaining part of this paper the ground field $k$ will be 
the algebraic closure of a finite field $\F_q$ with $q$ elements. In this
subsection we summarize some results from St\"ohr-Voloch's \cite[\S2]{sv}.

We keep the assumptions and the notations of the preceding subsection and 
we suppose that $X$ and $\cD$ are defined over $\F_q$. We let $\frx$ 
denote the Frobenius morphism (relative to $\fq$) on $X$. Then 
$X$ is equipped with {\it the $\fq$-Frobenius divisor associated to 
$\cD$}, namely
$$
S=S^{\cD}=\df(
{\rm det}
\begin{pmatrix} f_0\circ\frx   & \ldots & f_r\circ\frx   \\
                D^{\nu_0}_xf_0 & \ldots & D^{\nu_0}_xf_r \\
                  \vdots       & \vdots & \vdots         \\
            D^{\nu_{r-1}}_xf_0 & \ldots & D^{\nu_{r-1}f_r}
\end{pmatrix})+\sum_{i=0}^{r-1}\nu_i\df(dx)+(q+r)E\, ,
$$
where $x$ is 
a separating variable of $\fq(X)|\fq$, $f_0,\ldots,f_r$ is 
a $\fq$-base of $\cD'$, and 
$(\nu_0,\nu_1,\ldots,\nu_{r-1})$, called {\it the $\fq$-Frobenius orders 
of $\cD$}, is the minimum in the lexicographic order of the set 
$$
\{(j_0,\ldots,j_{r-1})\in \N^r: j_0<\ldots<j_{r-1};\ 
{\rm det}
\begin{pmatrix} f_0\circ\frx &\ldots & f_r\circ\frx\\
                D^{j_0}_xf_0 &\ldots & D^{j_0}_xf_r\\
                  \vdots     &\vdots & \vdots      \\
                D^{j_{r-1}}_xf_0 & \ldots & D^{j_{r-1}f_r}
\end{pmatrix}
\ \neq 0\}\, .
$$
One has $X(\fq)\subseteq \supp(S)$,
\begin{equation}\label{eq1.3}
\deg(S)=(\sum_{i=0}^{r-1}\nu_i)(2g-2)+(q+r)d\, ,
\end{equation}
$\nu_i=\epsilon_i$ for $i<I$, $\nu_i=\epsilon_{i+1}$  
for $i\ge I$, where $I=I^{(\cD,q)}$ is the smallest integer such that
$(f_0\circ\frx,\ldots,f_r\circ\frx)$ is a $\fq(X)$-linear combination of
the vectors $(D^{\epsilon_i}_xf_0,\ldots,D^{\epsilon_i}_xf_r)$ with
$i=0,\ldots,I$. For $P\in X(\fq)$ one also has
\begin{equation}\label{eq1.4}
\nu_i\le j_{i+1}(P)-j_1(P)\qquad i=0,\ldots,r-1\qquad \text{and}
\end{equation}
\begin{equation}\label{eq1.5}
v_P(S)\ge \sum_{i=0}^{r-1}(j_{i+1}(P)-\nu_i)\, .
\end{equation}
\subsection{A $\fq$-divisor from the Zeta Function}\label{1.3}
In this subsection we generalize \cite[Lemma 1.1]{fgt} and its
corollaries. Let $X\!\mid\!\fq$ be a curve and 
$h(t)=h_{X,q}(t)$ its {\it $h$-polynomial}, i.e. the characteristic 
polynomial of the Frobenius endomorphism $\frj$ of the Jacobian $\cJ$ 
(over $\fq$) of $X$. We let $\prod_{i=1}^{T}h^{r_i}_i(t)$ denote the 
factorization over $\Z [t]$ of $h(t)$. Because of the  
semisimplicity of $\frj$ and the faithfully of the representation 
of endomorphisms of $\cJ$ on the Tate module (cf. \cite[Thm.
2]{ta}, \cite[VI, \S3]{l}), we then have 
\begin{equation*}
\prod_{i=1}^{T}h_i(\frj)=0\quad \mbox{on}\ \cJ\, .\tag{$*$}
\end{equation*}
Throughout this subsection we set 
$$
\sum_{i=1}^{U }\alpha_i t^{U-i}+t^U:= 
\prod_{i=1}^{T}h_i(t)\, ,
$$
we assume that $X$ has at least one $\fq$-rational point, say $P_0$, and
set 
$$
\cD=\cD^{(X,q,P_0)}:= ||m|P_0|\qquad \text{with}\qquad 
m:=\prod_{i=1}^{T}h_i(1)\, .
$$
As $f\circ\frx=\frj\frx$, $f=f^{P_0}$ being the natural map from $X$ to 
$\cJ$ given by $P\mapsto [P-P_0]$, Eq $(*)$ is equivalent to 
\begin{equation}\label{eq1.6} 
\sum_{i=1}^{U}\alpha_i\frx^{U-i}(P)+\frx^U(P)\sim mP_0\, .
\end{equation}
This suggests the 
\smallskip

{\bf Problem.} Study the relation among the $\fq$-rational points, the 
Weierstrass points, the $\cD$-Weierstrass points, and the support of the
$\fq$-Frobenius divisor associated to $\cD$. 
\begin{lemma}\label{l1.1} 
\begin{enumerate} 
\item If $P, Q\in X(\fq)$, then $mP\sim mQ$. In particular, $|m|$ is a 
non-gap at each $P\in X(\fq)$ provided that $\# X(\fq)\ge 2$.
\item If ${\rm char}(\fq)\nmid m$ and $\#X(\fq)\ge 2g+3$, then there 
exists $P_1\in X(\fq)$ such that $|m|-1$ and $|m|$ are non-gap at $P_1$.
\end{enumerate}
\end{lemma}
\begin{proof} (1) It follows inmediately from (\ref{eq1.6}). (2) (The
proof is inspired by \cite[Prop.\! 1]{stix}.) By item (1) (and ${\rm
char }(\fq)\nmid m$) 
there exists a separable morphism 
$x:X \to \P^1(\bar \fq)$ with $\df(x)=|m|P_0-|m|Q$, $P_0, Q\in X(\fq)$. Let  
$n$ denote the number of unramified rational points for $x$. By the 
Riemann-Hurwitz formula we 
find that $n\ge \#X(\fq)-2g-2$ so that $n>0$ by the hypothesis
on $\#X(\fq)$. Thus there exists $\alpha \in \fq$, $P_1\in
X(\fq)\setminus \{P_0, Q\}$  
such that $\df(x-\alpha)=P_1+D-mQ$ with $P_1, Q \not\in \supp(D)$. Let
$y\in
\fq(X)$ be such that $\df(y)=|m|Q-|m|P_1$ (cf. item (1)). Then
from the rational function $(x-\alpha)y$ we obtain item (2).
\end{proof}
It follows that the definition of $\cD$ is independent of $P_0\in X(\fq)$ 
and the 
\begin{corollary}\label{cor1.1} 
\begin{enumerate}
\item If $\#X(\fq)\ge 2$, then $\cD$ is base-point-free.
\item If ${\rm char}(\fq)\nmid m$ and $\#X(\fq)\ge 2g+3$, then 
$\cD$ is simple. 
\end{enumerate}
\end{corollary}
Let us assume further that $\# X(\fq)\ge 2$ and that 
\begin{equation*}
m>0,\quad \alpha_1, \alpha_U\ge 1\quad\text{and}\quad 
\alpha_{i+1}\ge \alpha_i,\ \text{for}\ i=1,\ldots, U-1.\tag{Z}
\end{equation*}
\begin{remark*} There exist $h$-polynomials that do not satisfy $(Z)$; 
see, e.g. \cite{carb}.
\end{remark*}
We set $r:=\dim(\cD)$, i.e. $m_r(P)=m$ for each $P\in
X(\fq)$ (cf. Lemma \ref{l1.1}(1)). We keep the notations of the preceding 
subsections.
\begin{lemma}\label{l1.2} 
\begin{enumerate}
\item The $(\cD,P)$-orders for $P\in X(\fq)$ are $j_i(P)=m-m_{r-i}(P)$,
$i=0,\ldots,r$. 
\item $j_1(P)=1$ for $P\not\in X(\fq)$. In particular, $\epsilon_1=1$.
\item The numbers $1, \alpha_1,\ldots,\alpha_U$ are $\cD$-orders, so that
$r\ge U+1$.
\item If $\frx^i(P)\neq P$ for $i=1,\ldots, U+1$, then 
$\alpha_U$ is a non-gap at
$P$. In particular $\alpha_U$ is a generic non-gap of $X$.
\item If $\frx^i(P)\neq P$ for $i=1,\ldots, U$, but $\frx^{U+1}(P)=P$, 
then $\alpha_U-1$ is a non-gap at $P$.
\item If $g\ge \alpha_U$, then $X$ is non-classical.
\end{enumerate}
\end{lemma}
\begin{proof} Items (1), (2) and (3) can be proved as in \cite[Thm.
1.4(iii), Prop. 1.5(iii)]{fgt}. To prove (4), (5) and (6) we apply 
${(\frx)}_*$, as in \cite[IV, Ex.\! 2.6]{har}, to the equivalence in
(\ref{eq1.6}).
Then 
$$
\sum_{i=1}^{U}\alpha_i\frx^{U-i}(P)+\frx^U(P)\sim 
\alpha_1\frx^U(P)+\sum_{i=1}^{U-1}\alpha_{i+1}\frx^{U-i}(P) +\frx^{U+1}(P)
$$
so that
$$
\alpha_UP\sim
\frx^{U+1}(P)+(\alpha_1-1)\frx^U(P)+
\sum_{i=1}^{U-1}(\alpha_{i+1}-\alpha_i)\frx^{U-i}(P)\, .
$$
Now the results follow.
\end{proof}
\begin{corollary}\label{cor1.2} 
If $\# X(\fq)>q(m-\alpha_U)+1$, then $j_{r-1}(P)<\alpha_U$ for each $P\in
X(\fq)$.
\end{corollary}
\begin{proof} By Lemma \ref{l1.2}(1), $j_{r-1}(P)=m-m_1(P)$, and by  
\cite[Thm. 1(b)]{le}, $\# X(\fq)\le 1+qm_1(P)$. Then $j_{N-1}(P)\ge
\alpha_U$ implies $\# X(\fq)\le q(m-\alpha_U)+1$.
\end{proof}
\begin{corollary}\label{cor1.3} 
\begin{enumerate}
\item $\epsilon_r=\nu_{r-1}=\alpha_U$;
\item $X(\fq)\subseteq \supp(R^{\cD})$.
\end{enumerate}
\end{corollary}
\begin{proof} (1) By (\ref{eq1.1})(b), $\epsilon_{r-1}\le j_{r-1}$. Since
$\alpha_U$ is a $\cD$-order (Lemma \ref{l1.2}(3)) we conclude 
that $\epsilon_r=\alpha_U$. Now from (\ref{eq1.6}) and (\ref{eq1.2})
we have that $\frx(P)$ belong to the $(r-1)$-th osculating hyperplane;
thus $\nu_{r-1}=\epsilon_r$.

(2) By Lemma \ref{l1.2}(1) for each rational point $j_r(P)=m$. Since
$m>\alpha_U$, the result follows from (\ref{eq1.1})(c).
\end{proof}
\begin{corollary}\label{cor1.4} If $\# X(\fq)\ge q\alpha_U+1$, then
$m_1(P)=q$
for each $P\in X(\fq)$.
\end{corollary}
\begin{proof} Let $P\in X(\fq)$. By (\ref{eq1.1})(b), applied to the
canonical linear system, we have $m_1(P)\le m_1(Q)$, $Q$ being a generic
point of
$X$. Then $m_1(P)\le \alpha_U$ by Lemma \ref{l1.2}(4). On the other hand,
by \cite[Thm. 1(b)]{le} and the hypothesis on $\#X(\fq)$ we get $m_1(P)\ge
\alpha_U$ and we are done.
\end{proof}
\section{Maximal curves}\label{2}
In this section we shall be dealing with maximal curves over $\fq$ or
equivalently with curves over $\fq$ whose $h$-polynomial is 
$(t+\sqrt{q})^{2g}$, $g>0$ being the genus of the curve. We set 
$\l:=\sqrt{q}$. Then, by \S\ref{1.3}, each maximal curve $X\!\mid\!\fl$
is equipped with the linear system 
$$
\cD:=|(\ell+1)P_0| \qquad P_0\in X(\fl)\, ,
$$
which will be fixed throughout the entire section. Notice that $X$
satisfies the hypotheses in Lemma \ref{l1.1} and $(Z)$ in 
\S\ref{1.3}. By Corollary \ref{cor1.1}, $\cD$ is 
simple and base-point-free; by Lemma \ref{l1.2}(3), $\dim(\cD)\ge 2$; 
for each $P\in X$ relation (\ref{eq1.6}) reads (\cite[Corollary 1.2]{fgt})
\begin{equation}\label{eq2.1}
\l P+ \frx(P)\sim (\l+1)P_0\, .
\end{equation}
Then, for each $P\in X$
\begin{equation}\label{eq2.2}
m_0(P)=0<m_1(P)<\ldots <m_n(P)\le \l<m_{n+1}(P)\, ,
\end{equation}
where $n+1:=\dim(\cD)$. We keep the notations of the preceding section.
\subsection{Known results}\label{2.1} The results of this subsection have
been noticed in \cite[\S1]{fgt}.

2.1.1. By Corollary \ref{cor1.3}(1), $\epsilon_{n+1}=\nu_n=\l$; this 
together with \cite{ft} and \cite[Proof of Thm.\! 1]{hv} imply: 
$\nu_1=\l\Leftrightarrow n+1=2$, and $\nu_1=1\Leftrightarrow n+1\ge 3$.

2.1.2. Let $P\in X(\fl)$. By Lemma \ref{l1.2}(1),
$j_i(P)=\l+1-m_{n+1-i}(P)$, 
$i=0,1,\ldots,n+1$. Then $j_{n+1}(P)=m_{n+1}(P)=\l+1$. The case $i=\l$ in  
(\ref{eq1.5}) gives $j_1(P)=1$ so that $m_n(P)=\l$.

2.1.3. Let $P\not\in X(\fl)$. We assume $n+1\ge 3$ (the case 
$n+1=2$ has been studied in \cite{ft}; see also Theorem \ref{t2.1} here). 
>From (\ref{eq2.1}), $j_1(P)=1$, and   
$m_n(P)=\l$ whenever $\frx^2(P)\neq P$. Furthermore, by (\ref{eq1.1})  
and \cite[Thm.\! 1]{ho},  $m_{n-1}(P)=\l-1$ if $P$ is 
not a Weierstrass point of $X$.

Set $m_i:=m_i(P)$, $u_0:=1$; let $u=u(P), u_i=u_i(P)\in {\bar\fl}(X)$ such
that
$\df(u_i)=D_i-m_iP$, $P\not\in\supp(D_i)$, and
$\df(u)=\l P+\frx(P)-(\l+1)P_0$. Then
\begin{equation}\label{eq2.3}
(\l+1)P_0+\df(uu_i)=D_i+\frx(P)+(\l-m_i)P\, ,
\end{equation}
and so $0,1,\l-m_n,\ldots,\l-m_0$ are $(\cD,P)$-orders.

2.1.4. The $\cD$-orders. Let $\tilde m_i$ denote the $i$th
non-gap at a generic point. Then, by \S2.1.3, the 
$\cD$-orders are
\begin{align*}
& \epsilon_0=\l-\tilde m_n<\epsilon_1=1=\l-\tilde m_{n-1}<\ldots<
\epsilon_{n-J}=\l-\tilde m_J<\epsilon_{n-J+1}< \\
& \epsilon_{n-J+2}=\l-\tilde
m_{J-1}<\ldots<\epsilon_{n+1}=\l-\tilde m_0\, ,
\end{align*}
whit $J\in \N$, $1\le J\le n-1$.
\subsection{The $\fl$-Frobenius orders of $\cD$.}\label{2.2}
\begin{proposition*} With the notations of \S2.1.4, the 
$\fl$-Frobenius orders of $\cD$ are 
$$
\{\epsilon_0,\ldots,\epsilon_{n+1}\}\setminus\{\epsilon_{n-J+1}\}.
$$
\end{proposition*}
\begin{proof} For $P$ a generic point of $X$ let $u, u_i \in {\bar\fl}(X)$ 
be  as in \S2.1.3, and $v\in {\bar\fl}(X)$ such that
$$
(\l+1)P_0+\df(v)=\epsilon_{n-J+1}P+D_v\qquad P\not\in \supp(D_v)\, .
$$
>From this equation and (\ref{eq2.3}) we have that $
\{uu_n,uu_{n-1},\ldots,uu_J,v,uu_{J-1},\ldots,uu_0\}$ 
is a $(\cD,P)$-hermitian base of $\cD$. Hence, by (\ref{eq1.2}), $
\frx(P)\in L_{n-J+1}(P)$. Now the result follows from the
\begin{claim*}
\quad $\frx(P)\not\in L_{n-J}(P)$.
\end{claim*}
Indeed, if $\frx(P)\in L_{n-J}(P)$ then $\frx(P)\in
\supp(D_v)$; hence, by (\ref{eq2.1}), we would have 
$\l-\epsilon_{n-J+1}\in H(P)$-a contradiction.
\end{proof}
\begin{remark*} A slight modification of the above proof  
shows that each point $P\in \supp(S^\cD)\setminus X(\fl)$ is a Weierstrass 
point of $X$ (\cite{gt}).
\end{remark*}
\subsection{The morphism associated to $\cD$}\label{2.3} Let $\pi:X\to
\P^{n+1}(\bar\fl)$ be the morphism associated to $\cD$.
\begin{proposition*}\quad $\pi$ is a closed embedding, i.e. $X$ is
$\fl$-isomorphic to $\pi(X)$.
\end{proposition*}
\begin{proof}
By \cite[Prop.\! 1.9]{fgt}, we have to show that $m_n(P)=\l$ for each
$P\in X$. By \S2.1.2 we can assume that $P\not\in X(\fl)$. For such a $P$, 
suppose that $m_n(P)<\l$; then, as $j_1(P)=1$ and $j_{n+1}(P)=\l$, the
$(\cD,P)$-orders
would be
$$
0,1=\l-m_n(P),\l-m_{n-1}(P),\ldots,\l-m_1(P),\l\, .
$$
Hence (\ref{eq1.2}) and (\ref{eq2.3}) would imply $\frx(P)\in L_1(P)$. On
the other hand, the hyperplane corresponding to the function $uu_n(P)$ in 
\S2.1.3 contains $P$ (with multiplicity 1) and $\frx(P)$; thus it 
contains $L_1(P)$. This is a contradiction because the multiplicity of
$L_1(P)$ at $P$ is at least $j_2(P)=\l-m_{n-1}(P)\ge 2$. 
\end{proof}
Now \cite[Prop.\! 1.10]{fgt} can be state without the hypothesis on $\pi$:
\begin{corollary*} 
Let $X\!\mid\!\fl$ be a maximal curve of genus $g$ . For some $P\in
X(\fl)$ suppose that there exist $a, b\in H(P)$ such that all non-gaps
less than or equal to $\l+1$ are generated by $a$ and $b$. Then 
$H(P)=\langle
a,b\rangle$, so that $g=(a-1)(b-1)/2$.
\end{corollary*}
\subsection{Weierstrass points and maximal curves}\label{2.4}
In this section we show that each $\fl$-rational point of $X$ is a
Weierstrass point of the curve provided that $g$ is large enough. First we 
notice that (\ref{eq2.2}) implies $g\ge \l-n$ and that 
$$
g=\l-n\quad \Leftrightarrow\quad \{\l+1,\l+2,\ldots,\}\subseteq H(P)\quad
\forall P\in X\, .
$$
Since $\l$ is a non-gap for a non-Weierstrass point, cf. \S2.1.3, 
(\ref{eq2.2}) also implies (\cite[Prop. 1.7(i)]{fgt})
$$
\text{$X$ classical}\quad \Rightarrow\quad g=\l-n.
$$
We remark that $g=\l-n$ does not characterize classical maximal curves; 
see e.g. \cite[Prop. 1.8]{fgt}.

The following results are contained in the proof of \cite[Satz
II.2.5]{rai}.
\begin{lemma}\label{l2.1} Let $X\!\mid\!\fl$ be a maximal curve of genus
$g$ and $P$ a
non-Weierstrass point of $X$. If $\l+1\in H(P)$, then $\l+1,\ldots,2\l\in
H(P)$. In particular, since $\l\in H(P)$, $\{\l+1,\l+2,\ldots\}\subseteq
H(P)$ and $g=\l-n$.
\end{lemma}
\begin{proof} 
Let $i\in \{1,\ldots,\l\}$ such that $\l+i\not\in H(P)$; then
$\binom{\l+i-1}{\l} \not\equiv 0 \pmod{{\rm char}(\fl)}$. Hence, by the
$p$-adic criterion \cite[Corollary 1.9]{sv}, $\l+1\not\in H(P)$.
\end{proof}
\begin{corollary}\label{c2.1}\quad 
$
g=\l-n\quad \Leftrightarrow\quad \text{$\l+1$ is a non-gap at a
non-Weierstrass point of $X$}.$
\end{corollary}
\begin{corollary}\label{c2.2} If $g>\l-n$, then
$$
X(\fl)\subseteq \text{set of Weierstrass points of $X$}.
$$ 
\end{corollary}
\begin{proof}
It follows from the above corollary and \S2.1.2.
\end{proof}
\begin{remark*} There exists maximal curves with $g=\l-n$ where no
$\fl$-rational point is Weierstrass, see e.g. the remark after
Proposition \ref{p2.1}. The
hypothesis $g>\l-n$ is satisfied
if $g\ge \max{(2,\l-1)}$; indeed, $g=\l-1\le \l-n$ implies $n=1$, i.e.
$g=(\l-1)\l/2$ (\cite{ft}) and so $g\le 1$, a contradiction.
\end{remark*}
\begin{remark*} Let $X\!\mid\!\fl$ be non-hyperelliptic and maximal of 
genus $g$. Denote by $\cW$ the set of Weierstrass points of $X$ 
($=\supp(R^{K_X})$). Corollary \ref{c2.2} implies
$$
\# \cW >
\begin{cases}
16 & \text{if $g=3$},\\
25 & \text{if $g=6$},\\
\max{3(g+2), 4(g-1)} & \text{if $g\neq 3, 6$}.
\end{cases}
$$
Hence, we can use Pflaum's \cite[Corollary 2.6, Proof of Theorem 1.6]{pf}
to describe the isomorphism-class (over $\bar \fl$) and the automorphism
group ${\rm Aut}(X)$ (also over $\bar\fl$) of $X$ via Weierstrass points.
In fact, we conclude that the isomorphism-class of maximal curves is
determinated by their constellations of Weierstrass points and that 
$$
{\rm Aut}(X)\cong\{A\in PGL(g,\bar\fl): A\rho(\cW)=\rho(\cW)\},
$$
where $\rho:X\to \P^{g-1}(\bar\fl)$ is the canonical embedding. 
Notice that, as the morphism $\pi:X\to \P^{n+1}$ associated to $\cD$ is an 
embedding (\S\ref{2.3}), (\ref{eq2.1}) implies
$$
{\rm Aut}_{\fl}(X)\cong \{A\in PGL(n+1,\fl): A\pi(X)=\pi(X)\}\, .
$$
\end{remark*}
\subsection{On the genus of maximal curves}\label{2.5} It has been
noticed in  \cite{ft} that the genus $g$ of a maximal curve
$X\!\mid\!\fl$ satisfies
$$
g\le (\l-1)^2/4\qquad\text{or}\qquad g=(\l-1)\l/2\, ,
$$
which was conjectured by Stichtenoth and Xing (cf. \cite{stix}).  
Moreover, we have the
\begin{theorem}\label{t2.1}
The following statements are equivalent
\begin{enumerate}
\item $X$ is $\fl$-isomorphic to $y^\l+y=x^{\l+1}$ (the Hermitian curve
over
$\fl$); 
\item $X\!\mid\!\fl$ maximal with $g>(\l-1)^2/4$;
\item $X\!\mid\!\fl$ maximal with ${\rm dim}(\cD)=2$.
\end{enumerate}
\end{theorem}
It is well known that $y^\l+y=x^{\l+1}$ is a maximal curve over $\fl$ of 
genus $(\l-1)\l/2$; $(2)\Rightarrow (3)$ follows by Castelnuovo's genus
bound for curves in projective spaces \cite[Claim 1]{ft}; $(3)\Rightarrow
(1)$ is the main result of \cite{ft}. Next we write a new proof of this
implication.
\begin{proof} 
$(3)\!\Rightarrow\! (1):$ By \S\ref{2.1},  
$(\epsilon_0,\epsilon_1,\epsilon_2)=(0,1,\l)$, $(\nu_0,\nu_1)=(0,\l)$, and 
$(j_0(P),j_1(P),j_2(P))=(0,1,\l+1)$ for each $P\in X(\fl)$. Hence 
(\ref{eq1.1})(a)(c) imply $g=(\l-1)\l/2$. Let $x, y\in \fl(X)$ with
$\dinf(x)=\l P_0$ and 
$\dinf(y)=(\l+1)P_0$. Then 
$H(P_0)=\langle \l,\l+1\rangle$ and so $\df(dx)=(2g-2)P_0\ (*)$, 
because $H(P_0)$ is symmetric. By $\nu_1=\l$ we have an equation of type
(cf. \S\ref{1.2})
\begin{equation}\label{eq2.4}
y^{\l^2}-y=f(x^{\l^2}-x)\ ,
\end{equation}
with $f:=D^{1}y$ (derivation with respect to $x$); by $\epsilon_2=\l$ we 
have the following two-rank matrices (cf. \S\ref{1.1})
$$
\left( \begin{array}{ccc}
1 & x   & y\\
0 & 1   & D^{(1)}y\\
0 & 0   & D^{(j)}y
\end{array} \right),
\qquad 2\le j<\epsilon_2=\l\, .
$$
By $(*)$ and (\ref{eq2.4}), $\dinf(f)=\l^2P_0$. Now 
$D^{(j)}y=0$ for $2\le j<\l$ and (\ref{eq2.4}) imply 
$D^{(j)}f=0$ for $1\le j<\l$. Thus, by  
\cite[Satz 10]{hasse}, there exists $f_1\in \fl(X)$ such that  
$f=f_1^{\l}$. Therefore $f_1=ax+b$ with $a, b\in
\fl$, $a\neq0$, and after some $\fq$-linear transformations we obtain an
equation of type
$$
y_1^{\l}+y_1-x_1^{\l+1}=(y_1^\l+y_1-x_1^{\l+1})^{\l}\, ,
$$
with $x_1, y_1\in \fq(X)$. Now the proof follows.
\end{proof}
\begin{remark*} Let $X$ be the Hermitian curve over $\fl$. From the above
proof we have $\cK_X=(\l-2)\cD$ and $(j_0(P),j_1(P),j_2(P))=(0,1,\l)$ for
each $P\not\in X(\fl)$. Then the $(\cK_X,P)$-orders contains 
$ \{a+b\l: a,b\ge 0, a+b\le \l-2\}$ (resp. 
$\{a+b(\l+1): a,b\ge 0, a+b\le \l-2\}$) if $P\not\in X(\fl)$ (resp. $P\in
X(\fl)$). Since these sets have cardinality equal to $g=(\l-1)\l/2$, these 
are the $(\cK_X,P)$-orders; hence 
$$
X(\fl)=\supp(R^\cD)=\supp(R^{\cK_X} (= \text{set of Weierstrass points of
$X$})\, ,
$$
and we have another proof of the fact that $X$ is non-classical (compare
with Lemma \ref{l1.2}(6)). 
The above computations have been carried out in \cite{gv}. We mention  
that the first examples of non-classical curves were obtained from
certain Hermitian 
curves (see \cite{sch}).
\end{remark*}
Now let us consider the following property for the maximal curve 
$X\!\mid\!\fl$ of genus $g$ (recall that $n+1=\dim(\cD)$):
\begin{equation}\label{eq2.5}
\exists P_1\in X(\fl)\quad \exists m\in H(P_1)\quad \text{such that $mn\le
\l+1$}\, .
\end{equation}
Then $mn=\l+1$ or $mn=\l$. In both cases the hypothesis of the corollary 
in \S\ref{2.3} is satisfied; in particular, $g=(\l-1)(m-1)/2$ or
$g=\l(m-1)/2$. 

{\bf Case $mn=\l+1$.} This occurs iff $X$ is $\fl$-isomorphic to 
$y^\l+y=x^{(\l+1)/n}$ \cite[Thm. 2.3]{fgt}) (so that $g= 
\frac{\l-1}{2}(\frac{\l+1}{n}-1)$). Hence there exists maximal curves of
genus $(\l-1)^2/4$ and  
indeed, $y^\l+y=x^{(\l+1)/2}$ is the unique maximal curve (up to
$\fl$-isomorphism) having such a genus (\cite[Thm.
3.1]{fgt}). Before we consider the case $mn=\l$ we 
prove an analogue of $(1)\Leftrightarrow (2)$ of Theorem \ref{t2.1}.
\begin{proposition}\label{p2.1} Let $X\!\mid\!\fl$ be a maximal curve of
genus $g$ and assume that $\l$ is odd. Then the following statements are
equivalent
\begin{enumerate}
\item $X$ is $\fl$-isomorphic to $y^\l+y=x^{(\l+1)/2}$;
\item $(\l-1)(\l-2)/4<g\le(\l-1)^2/4$.
\end{enumerate}
Item (1) (or (2)) implies\quad $3.\, \, \dim(\cD)=3$.
\end{proposition}
\begin{proof} 
We already noticed that $(1)\Rightarrow (2)$ and $(1)\Rightarrow
(3)$.
That $(2)\Rightarrow (3)$ follows by Castelnuovo's genus bound for curves
in projective spaces \cite{c}, \cite[p.\! 116]{acgh}, \cite[Corollary
2.8]{ra}.

$(2)\Rightarrow (1):$ The cases $\l\le 5$ are trivial, so let $\l>5$.  
According to \cite[Thm.\! 2.3]{fgt}, we look for a rational point $P$
such that there exists $m\in H(P)$ with $2m=\l+1$. Let 
$m_1:=m_1(P)<\l<\l+1$ be the first three positive non-gaps at $P\in
X(\fl)$. By
\S2.1.2 the $(\cD,P)$-orders are $0,1,j=\l+1-m_1,\l+1$. Notice that $\l$
odd implies $2m_1\ge \l+1$ and hence that $2j\le \l+1$. 

Set $2\cD:=|2(\l+1)P_0|$; $\dim(\cD)=3$ implies $\dim(2\cD)\ge
8$; the lower bound on $g$ implies (once again via Castelnuovo's bound) 
$\dim(2\cD)=8$. The $(2\cD,P)$-orders ($P\in X(\fl)$) contains the set
$\{0,1,2,j,j+1,2j,\l+1,\l+2,\l+j+1,2\l+2\}$; 
therefore $\dim(2\cD)=8$ implies $j=2$ (i.e. $m_1(P)=\l-1$) or $2j=\l+1$
(i.e.
$m_1(P)=(\l+1)/2$) and we have to show that it is not possible to
have $m_1(P)=\l-1$ for each $P\in X(\fl)$. 

Suppose that $m_1(P)=\l-1$ for each $P\in X(\fl)$. Then the 
$\cD$-orders are $0,1,2,\l$ and so $v_P(R_1)=1$ for each $P\in X(\fl)$ 
($R_1$ being the $\cD$-ramification divisor). Then, by (\ref{eq1.1}),
$$
\deg(R_1)-\#X(\fl)=3(2g-2)-(\l-3)(\l+1)\, .
$$
\begin{claim*}
\quad For each $P\in \supp(R_1)\setminus X(\fl)$ the $(\cD,P)$-orders are
$0,1,(\l+1)/2,\l$.
\end{claim*}
In fact, for such a $P$ the $(\cD,P)$-orders
are $0,1,i,\l$ with $i=i(P)>2$, and 
$$
\{0,1,2,i,i+1,2i,\l,\l+1,\l+i,2\l\}
$$
is contained in the $(2\cD,P)$-orders; thus $\dim(2\cD)=8$ implies
$i\in\{(\l+1)/2,\l-1\}$. 

Suppose that $i=\l-1$; by (\ref{eq2.1}) there exists $Q_1, Q_2\in
X\setminus\{P\}$ such that $P+\frx(P)\sim Q_1+Q_2$, i.e. $X$ is
hyperelliptic. This implies $g\le (\l-1)/2$ (see e.g. \cite[Thm. 
1(b)]{le}) and
from the hypothesis on $g$ that $\l< 4$, a contradiction.

By the claim and (\ref{eq1.1})(d), for each $P\in
\supp(R_1)\setminus X(\fl)$, $v_P(R_1)=(\l-3)/2$ and
$$
A:=\#(\supp(R_1)\setminus X(\fl))=\frac{6(2g-2)}{\l-3}-2(\l+1)\, .
$$
With the above computations we analize $(2\cD,P)$-orders for $P\in
\supp(R_1)$. We have:
$$
\text{$(2\cD,P)$-orders}\ =
\begin{cases}
0,1,2,3,4,\l+1,\l+2,\l+3,2\l+2 & \text{if $P\in X(\fl)$};\\
0,1,2,(\l+1)/2,(\l+3)/2,\l,\l+1,(3\l+1)/2,2\l  & \text{if $P\not\in 
X(\fl)$}.
\end{cases}
$$
Denote by $R_2$ the $2\cD$-ramification divisor. Being 
$0,1,2,3,4,\l,\l+1,\l+2,2\l$ the $2\cD$-orders we then have 
$$
v_P(R_2)\ge
\begin{cases}
5 & \text{if $P\in X(\fl)$};\\
(3\l-13)/2 & \text{if $P\in \supp(R_1)\setminus X(\fl)$}.
\end{cases}
$$
Then, again by (\ref{eq1.1}), 
$$
\deg(R_2)=(5\l+13)(2g-2)+18(\l+1)\ge 5\# X(\fl)+\frac{3\l-13}{2} A\, ,
$$
which implies $2g-2\ge (\l-3)(\l+1)/2$, i.e. 
$2g-2=(\l-3)(\l+1)/2$ due to the upper bound on $g$. By \cite[Thm.\! 
3.1]{fgt} we then conclude that $A=0$, i.e. $2g-2=(\l-3)(\l+1)/3$, a
contradiction. This complete the proof.
\end{proof}
\begin{remark*} Let $X\!\mid\!\fl$ be maximal of genus $g$ and suppose 
that $(\l-1)(\l-2)/6<g\ge (\l-1)^2/4$. If $X \not\cong 
y^\l+y=x^{(\l+1)/2}$, then $(\l-1)(\l-2)/6<g\le (\l-1)(\l-2)/4$ by the 
last proposition. Cossidente and Korchmaros \cite{ck} constructed a 
maximal curve $X\!\mid\!\fl$ with $g=(\l+1)(\l-2)/6$ and $\l\equiv
2\pmod{3}$. By Castelnuovo's 
genus bound, the linears system $\cD$ of this curve satisfies 
$\dim(\cD)=3$; so this example shows that $(3)$ does not imply $(1)$. 
\end{remark*}
\begin{remark*} For $\l\equiv 3\pmod{4}$ there are at least two non 
$\bar\fl$-isomorphism maximal curves over $\fl$ with
$g=\frac{\l-1}{2}(\frac{\l+1}{4}-1)$ and $\dim(\cD)=5$, namely
$$
(1)\ \ y^\l+y=x^{(\l+1)/4}\qquad\text{and}\qquad (2)\ \ 
x^{(\l+1)/2}+y^{(\l+1)/2}=1.
$$
Indeed, the curve in (2) admits points $P$ with $H(P)=\langle
(\l-1)/2,(\l+1)/2\rangle$ (e.g. $P$ over a root of $x^{\l+1}=1$) and is 
well known that such semigroups cannot be realized by (1) (\cite{gv}).

These examples show that one cannot expect the uniqueness of a maximal
curve just by means of a given genus. It also shows that the hypothesis on
non-gaps of \cite[Prop.\! 1.10, Thm.\! 2.3]{fgt} cannot be relaxed.

The curve (2) have been considered by Hirschfeld and Korchmaros \cite{hk}. 
They noticed an interesting bound for the number of rational points
of a curve; the curve in (2) attains such a bound.    
\end{remark*}
\begin{remark*} In view of the above examples and \cite{ck}'s letter is 
reasonable to make the 
following conjectures. Let $X\!\mid\!\fl$ be a maximal curve of genus
$g$.  

(1) Let $\l$ be odd. If $\l\not\equiv 2\pmod{3}$, then
$(\l-1)(\l-2)/6 <g\le (\l-1)^2/4$ iff 
$g=(\l-1)^2/4$ (i.e. $X$ is $\fl$-isomorphic to $y^\l+y=x^{(\l+1)/2}$). If
$\l\equiv 2\pmod{3}$, then $(\l-1)(\l-2)/6<g\le (\l-1)^2$ iff
$g=(\l+1)(\l-2)/6$ or $g=(\l-1)^2/4$.

(2) With the exception of finitely many $\l$'s and if $\l\equiv 2
\pmod{3}$, then $(\l-1)(\l-3)/8<g\le (\l-1)(\l-2)/6$
iff $X$ is $\fl$-isomorphic to $y^\l+y=x^{(\l+1)/3}$ (in particular
$g=(\l-1)(\l-2)/6$).
\end{remark*}
{\bf Case $mn=\l$.} Now we assume (\ref{eq2.5}) with $mn=\l$. To
begin we notice that the quotient of the Hermitian curve by a certain
automorphism has 
a plane model of type $F(y)=x^{\l+1}$, $F$ being an additive polynomial.
These curves provide examples of maximal curves for this case. It has been 
conjectured in \cite{fgt} that $X$ 
is $\fl$-isomorphic to the above plane model with
$\deg(F)=m$; the fact that $g=\l(m-1)/2$ may provide evidence for
this conjecture. Next we state another proof of this fact, where 
is implicitely outlined a method to find a plane model for $X$:

By Theorem \ref{t2.1} we can assume 
$n>1$. Let $x,y\in \fl(X)$ such that $\dinf(x)=m$ and
$\dinf(y)=\l+1$. 
\begin{claim*} For each $\alpha\in \fl$, $\# x^{-1}(\alpha)=m$ and
$x^{-1}(\alpha)\subseteq X(\fl)$.
\end{claim*}
This implies $g=\l(m-1)/2$ because $\deg(\df(x^{\l^2}-x))=0$ gives $
\l^2+2\l g=\l^2m$.\newline 
The claim follows from two facts: 
\begin{fact}\label{f1} For each $P\neq P_1$, $\#x^{-1}(x(P))=m$.
\end{fact}
\begin{proof} {\it (Fact \ref{f1})} Let $P\neq P_1$ and for $x(P)=\alpha
\in
\bar\fl$ set
$e=v_P(x-\alpha)$. We have to show that $e=1$. Writing 
$\df(x-\alpha)=eP+D_P-mP_1$ with $P\not\in\supp(D_P)$, we then see
that $e,\ldots,en$ are $(\cD,P)$-orders. If $e>1$, then $j_{n+1}(P)=en$
because 1 is a $(\cD,P)$-orders (cf. \S\ref{2.1}). This implies $P\not\in
X(\fl)$ and so $e=m$ because $\l=mn$. Consequently $mP\sim mP_1$ so that 
$\l P\sim\l P_1$. Then by (\ref{eq2.1}) we get $\frx(P)\sim P_1$, a
contradiction because $g>0$. This finish the proof of Fact 1.
\end{proof}
Let $P\neq P_1$. From the above proof, the 
$(\cD,P)$-orders are 
$0,1,\ldots,n,\l+1$ (resp. $0,1,\ldots,n,\l$) if $P\in X(\fl)$ (resp.
$P\not\in X(\fl)$). Hence the $\cD$-orders are $0,1,\ldots,n,\l$ so that
$\supp(R^{\cD})=X(\fl)$ with $v_P(R^{\cD})=1$, $P\in
X(\fl)\setminus\{P_1\}$ , cf. (\ref{eq1.1}). 

Now the morphism $\pi$ associated to $\cD$ can be defined by
$(1:x:\ldots:x^n:y)$; so $v_P(D^\l y)=1$ for $P\in
X(\fl)\setminus P_1$ and $v_P(D^\l y)=0$ for $P\not\in X(\fl)$ (derivation
with respect to $x$); cf. \S\ref{1.1}.
\begin{fact}\label{f2} There exists $f\in \fl(X)$ regular outside $P_1$
such that $D^\l y= f(x^{\l^2}-x)$.
\end{fact}
\begin{proof} {\it (Fact \ref{f2})} By (\ref{eq1.4}) the $\fl$-Frobenius
orders
are $0,1,\ldots,n-1,\l$. Hence (cf. \S\ref{1.2})
$$
{\rm det}
\begin{pmatrix} 1 & x^{\l^2} & \ldots & x^{\l^2 n} &  y^{\l^2}\\
                1 &    x     & \ldots &    x^n     &  y       \\
                0 &   D^1x   & \ldots &   D^1x^n   &  D^1y    \\
          \vdots  & \vdots   & \vdots & \vdots     & \vdots   \\
                0 &   D^nx   & \ldots &   D^nx^n   &  D^n  
\end{pmatrix} =0
$$
and so
$$
y-y^{\l^2}=\sum_{i=1}^{n}(x^i-x^{i\l^2 })h_i
$$
with $h_i\in \fl(X)$ regular in $X\setminus\{P_1\}$ for each $i$. Fact
\ref{f2} 
now follows by applying $D^\l$ to the above equation and using the
following properties:
\begin{itemize}
\item For $r$ a power of a prime we have $D^j f^r=(D^{j/r}f)^r$ if $r\mid
j$ and $D^j f^r=0$ otherwise;
\item the fact that $n$ is a power of a prime implies $D^i x^n=0$ for
$i=1,\ldots,n-1$;
\item $D^i\circ D^j=\binom{i+j}{i}D^{i+j}$\, .
\end{itemize}
\end{proof}
\section{On the Deligne-Lusztig curve associated to the Suzuki
group}\label{3}
In this section we prove Theorem \ref{B} stated in the introduction. 
Throughout
we let $q_0:=2^s$, $q:=2q_0$, $X\!\mid\!\fq$ a curve of genus $g$ with
$$
g=q_0(q-1)\qquad \text{and}\qquad \# X(\fq)=q^2+1\, .
$$ 
Then, by Serre-Weil's explicit formulae (cf. \cite{se}, \cite{h}), the
$h$-polynomial of $X$ is $(t^2+2q_0t+q)^g$. Hence, by \S\ref{1.3}, $X$ is
equipped with the base-point-free simple linear 
system $\cD:=|(q+2q_0+1)P_0|$, with $P_0\in X(\fq)$. Here for each $P\in
X$, (\ref{eq1.6}) reads
\begin{equation}\label{eq3.1} 
qP+2q_0\frx(P)+\frx^2(P)\sim (q+2q_0+1)P_0\, .
\end{equation} 
We notice that $X$ satisfies the hypothesis in Lemma 
\ref{l1.1} and $(Z)$ in \S\ref{1.3}. As a first consequence we have the
\begin{lemma}\label{l3.1} Let $X$ be as above. Then for each $P\in
X(\fq)$, $m_1(P)=q$.  
\end{lemma}
Next we apply \cite{sv} to $\cD$; we keep the notation in \S\ref{1} and 
set $r:=\dim(\cD)$. The key property of $\cD$ 
will be the fact that $\frx(P)$ belongs to the tangent 
line for $P$ generic (Lemma \ref{l3.4}(1)). For $P\in X(\fq)$, Lemma
\ref{l1.2}(1) gives 
$m_r(P)=q+2q_0+1$ and 
\begin{equation}\label{eq3.2}
j_i(P)=m_r(P)-m_{r-i}(P)\quad \mbox{for}\ \ i=0,\ldots,r\, .  
\end{equation}
This together with Lemma \ref{l3.1} imply 
\begin{equation}\label{eq3.3}
j_r=(P)=q+2q_0+1,\quad j_{r-1}(P)=2q_0+1\qquad  (P\in
X(\fq))\, .
\end{equation} 
By Lemma \ref{l1.2}(3) $1, 2q_0, q$ are $\cD$-orders, so $r\ge 3$ and 
$\epsilon_1=1$. By Corollary \ref{cor1.3}(1) and (\ref{eq1.1})(b),
\begin{equation}\label{eq3.4}
\epsilon_r=\nu_{r-1}=q\qquad\text{and}\qquad 2q_0\le\epsilon_{r-1}\le
1+2q_0\, .
\end{equation} 
\begin{lemma}\label{l3.2}\quad $\epsilon_{r-1}=2q_0$.  
\end{lemma}
\begin{proof} Suppose that $\epsilon_{r-1}> 
2q_0$. Then, by (\ref{eq3.4}), $\epsilon_{r-2}=2q_0$ and 
$\epsilon_{r-1}=2q_0+1$. By (\ref{eq1.4}) and (\ref{eq3.3}), 
$\nu_{r-2}\le 2q_0=\epsilon_{r-2}$. Thus the $\fq$-Frobenius orders of 
$\cD$ would be $\epsilon_0,\epsilon_1,\ldots,\epsilon_{r-2}$,  
and $\epsilon_r$. By (\ref{eq1.5}) and (\ref{eq1.4}), for each $P\in
X(\fq)$
\begin{equation}\label{eq3.5}
v_P(S)\ge \sum_{i=1}^{r}(j_i(P)-\nu_{i-1})\ge
(r-1)j_1(P)+1+2q_0\ge r+2q_0\ ,
\end{equation} 
Thus $\deg(S)\ge (r+2q_0)\#X(\fq)$, and
from (\ref{eq1.3}), the identities $2g-2=(2q_0-2)(q+2q_0+1)$ and 
$\#X(\fq)=(q-2q_0+1)(q+2q_0+1)$ we obtain 
$$
\sum_{i=1}^{r-2}\nu_i=\sum_{i=1}^{r-2}\epsilon_i\ge (r-1)q_0\, .  
$$ 
Now, as $\epsilon_i+\epsilon_j\le \epsilon_{i+j}$ for $i+j\le r$ 
\cite[Thm. 1]{e}, then we would have 
$$ 
(r-1)\epsilon_{r-2}\ge 2\sum_{i=0}^{r-2}\epsilon_i\, ,  
$$ 
and hence $\epsilon_i+\epsilon_{r-2-i}=\epsilon_{r-2}$ for
$i=0,\ldots,r-2$. In particular,
$\epsilon_{r-3}=2q_0-1$ and by the $p$-adic criterion (cf. \cite[Corollary
1.9]{sv} we would have $\epsilon_i=i$ for $i=0,1,\ldots,r-3$. These facts
imply $r=2q_0+2$. Finally, we are going to see that this is a
contradiction via Castelnuovo's genus bound \cite{c}, \cite[p.\! 116]{acgh},
\cite[Corollary 2.8]{ra}. Castelnuovo's formula applied to $\cD$ implies
$$
2g=2q_0(q-1)\le \frac{(q+2q_0-(r-1)/2)^2}{r-1}\, .
$$
For $r=2q_0+2$ this gives $2q_0(q-1)< (q+q_0)^2/2q_0=q_0q+q/2+q_0/2$, a
contradiction.
\end{proof}
\begin{lemma}\label{l3.3}
There exists $P_1\in X(\fq)$ such that 
$$
\left\{ \begin{array}{ll}
j_1(P_1)=1         & {} \\
j_i(P_1)=\nu_{i-1}+1 & \mbox{if}\ i=2,\ldots, r-1.
\end{array}\right.
$$
\end{lemma}
\begin{proof} By (\ref{eq3.5}), it is enough to show that there exists
$P_1\in X(\fq)$ such that $v_{P_1}(S)=r+2q_0$. Suppose that $v_P(S)\ge
r+2q_0+1$ for each 
$P\in X(\fq)$. Then by (\ref{eq1.3}) we would have that
$$
\sum_{i=0}^{r-1}\nu_i \ge q+rq_0+1\ ,
$$
and, as $\epsilon_1=1$, $\nu_{r-1}=q$ and $\nu_i\le \epsilon_{i+1}$, that 
$$
\sum_{i=0}^{r-1}\epsilon_i \ge rq_0+2\, .
$$
By \cite[Thm. 1]{e} (or \cite{ho}), we then would conclude that 
$r\epsilon_{r-1}\ge 2rq_0+4$, i.e. $\epsilon_{r-1}>2rq_0$, a contradiction 
with the previous lemma.
\end{proof}
\begin{lemma}\label{l3.4}
\begin{enumerate} 
\item \ $\nu_1>\epsilon_1=1$.
\item \ $\epsilon_2$ is a power of two.
\end{enumerate}
\end{lemma}
\begin{proof} Statement (2) is consequence of the $p$-adic criterion 
\cite[Corollary 1.9]{sv}.  
Suppose that $\nu_1=1$. Then by Lemma \ref{l3.1}, Lemma \ref{l3.3}, 
(\ref{eq3.3}) and (\ref{eq3.2}) there would be a point $P_1\in X$ such
that $H(P_1)$ would contain the semigroup $H:=\langle q, q+2q_0-1, q+2q_0, 
q+2q_0+1\rangle $. Then $g\le  \#(\N\setminus
H)$, a contradiction as follows from the remark below. 
\end{proof}
\begin{remark*} Let $H$ be the semigroup defined above. We are going to
show that $\tilde g:= \#(\N \setminus H)=g-q_0^2/4$. 
To begin, notice that  
$L:=\cup_{i=1}^{2q_0-1}L_i$ is a complete system of residues
module $q$, where
$$
\begin{array}{lll}
L_i & = &
\{iq+i(2q_0-1)+j: j=0,\ldots,2i\}\quad  \mbox{if}\ \ 1\le i\le q_0-1,\\
L_{q_0} & = & \{q_0q+q-q_0+j:j=0,\ldots,q_0-1\},\\
L_{q_0+1} & = & \{(q_0+1)q+1+j:j=0,\ldots,q_0-1\},\\ 
L_{q_0+i} & = & 
\{(q_0+i)q+(2i-3)q_0+i-1+j: 
j=0,\ldots,q_0-2i+1\}\cup\\
          &   & \{(q_0+i)q+(2i-2)q_0+i+j: j=0,\ldots q_0-1\}\quad  
\mbox {if}\ \ 2\le i\le q_0/2,\\
L_{3q_0/2+i} & = & 
\{(3q_0/2+i)q+(q_0/2+i-1)(2q_0-1)+q_0+2i-1+j:\\
             &  &  
j=0,\ldots,q_0-2i-1\}\quad \mbox {if}\ \ 1\le i\le q_0/2-1.
\end{array}
$$
Moreover, for each $\l \in L$, $\l \in H$ and $\l-q\not\in H$. Hence
$\tilde g$ can be computed by summing up the coefficients of $q$ from the 
above list (see e.g. \cite[Thm. p.3]{sel}), i.e.
$$
\begin{array}{lll}
\tilde g & = & \sum_{i=1}^{q_0-1}i(2i+1)+q_0^2+(q_0+1)q_0+
\sum_{i=2}^{q_0/2}(q_0+i)(2q_0-2i+2)+\\
         &   & 
\sum_{i=1}^{q_0/2-1}(3q_0/2+i)(q_0-2i)=q_0(q-1)-q_0^2/4\, .
\end{array}
$$
\end{remark*}
In the remaining part of this section let $P_1\in X(\fq)$ be a point  
satisfying Lemma \ref{l3.3}; we set $m_i:= m_i(P_1)$ and denote by $v$   
the valuation at $P_1$. 

The item (1) of the last lemma implies $\nu_i=\epsilon_{i+1}$ for 
$i=1,\ldots, r-1$. Therefore from (\ref{eq3.2}),  
(\ref{eq3.3}) and Lemma \ref{l3.3},
\begin{equation}\label{eq3.6}
\left\{\begin{array}{ll}
m_i=2q_0+q-\epsilon_{r-i} & \mbox{if}\ i=1,\ldots r-2\\
m_{r-1}=2q_0+q,\ \ m_r=1+2q_0+q. & {}
\end{array}\right.
\end{equation}
Let $x, y_2,\ldots, y_r\in \fq(X)$ be such that $\dinf(x)=m_1P_1$, and 
$\dinf (y_i)=m_i P_1$ for $i=2,\ldots, r$. The fact that $\nu_1>1$ means
that the following matrix has rank two (cf. \S\ref{1.2})
$$
\left( \begin{array}{ccccc}
1 & x^q & y_2^q &\ldots &y_r^q\\
1 & x   & y_2   &\ldots &y_r\\
0 & 1   & D^{(1)}y_2   &\ldots& D^{(1)}y_r
\end{array} \right)\, .
$$
In particular, 
\begin{equation}\label{eq3.7}
y_i^q-y_i= D^{(1)}y_i(x^q-x) \quad \text{for}\ \ i=2,\ldots, r.
\end{equation}
\begin{lemma}\label{l3.5} 
\begin{enumerate}
\item For $P\in X(\fq)$, the divisor $(2g-2)P$ 
is canonical, i.e. the Weierstrass semigroup at $P$ is symmetric. 
\item  Let $m\in H(P_1)$. If $m<q+2q_0$, then $m\le
q+q_0$. 
\item For $i=2,\ldots,r$ there exists $g_i\in \fq(X)$ such that $
D^{(1)}y_i=g_i^{\epsilon_2}$. 
Furthermore, $\dinf(g_i)=\frac{qm_i-q^2}{\epsilon_2}P_1$.
\end{enumerate}
\end{lemma}
\begin{proof} (1) Since $2g-2=(2q_0-2)(q+2q_0+1)$, by (\ref{eq3.1}) 
we can assume $P=P_1$.  Now the case $i=r$ of 
Eqs. (\ref{eq3.7}) implies $v(dx)=2g-2$ and we are done.

(2) By (\ref{eq3.6}), $q, q+2q_0$ and $q+2q_0+1\in H(P_1)$. Then the
numbers
$$
(2q_0-2)q+q-4q_0+j\qquad j=0,\ldots,q_0-2
$$
are also non-gaps at $P_1$. Therefore, by the symmetry of $H(P_1)$,
$$
q+q_0+1+j\qquad j=0,\ldots,q_0-2
$$
are gaps at $P_1$ and the proof follows.

(3) Set $f_i:= D^{(1)}y_i$. By the product rule applied to
(\ref{eq3.7}),\\  
$D^{(j)}y_i=(x^q-x)D^{(j)}f_i+D^{(j-1)}f_i$ for $1\le j<q$. 
Then, $D^{(j)}f_i=0$ for $1\le 
j<\epsilon_2$, because the matrices
$$
\left( \begin{array}{ccccc}
1 & x   & y_2   &\ldots &y_r\\
0 & 1   & D^{(1)}y_2   &\ldots& D^{(1)}y_r\\
0 & 0   & D^{(j)}y_2   &\ldots& D^{(j)}y_r
\end{array} \right),
\quad 2\le j<\epsilon_2
$$
have rank two (cf. \S\ref{1.1}). Consequently, as $\epsilon_2$ is a power
of 2 (Lemma \ref{l3.4}(2)), by \cite[Satz 10]{hasse}, 
$f_i=g_i^{\epsilon_2}$ for some $g_i\in \fq(X)$. Finally, from the proof
of item (1) we have that $x-x(P)$ is a local 
parameter at $P$ if $P\neq P_1$. Then, by the election of the $y_i$'s, 
$g_i$ has no pole but in $P_1$, and from (\ref{eq3.7}),   
$v(g_i)=-(qm_i-q^2)/\epsilon_2$. 
\end{proof}
\begin{lemma}\label{l3.6}\quad $r=4$ and $\epsilon_2=q_0$.
\end{lemma}
\begin{proof} 
We know  that $r\ge 3$; we claim that $r\ge 4$; in 
fact, if $r=3$ we would have $\epsilon_2=2q_0$, $m_1=q$, $m_2=q+2q_0$,
$m_3=q+2q_0+1$, and hence $v(g_2)=-q$ ($g_2$ being as in Lemma  
\ref{l3.5}(3)). Therefore, after some $\fq$-linear transformations, the
case $i=2$ of (\ref{eq3.7}) reads 
$$
y_2^q-y_2=x^{2q_0}(x^q-x)\, .
$$
Now the function $z:= y_2^{q_0}-x^{q_0+1}$ satisfies
$z^q-z=x^{q_0}(x^q-x)$ and we find that $q_0+q$ is
a non-gap at $P_1$ (cf. \cite[Lemma 1.8]{hsti}). This contradiction
eliminates the case $r=3$.

Let $r\ge 4$ and $2\le i\le r$. By Lemma \ref{l3.5}(3)  
$(qm_i-q^2)/\epsilon_2\in H(P_1)$, and since $(qm_i-q^2)/\epsilon_2\ge
m_{i-1}\ge q$, by (\ref{eq3.6}) we have 
$$
2q_0\ge \epsilon_2 +\epsilon_{r-i}\qquad \mbox{for}\ i=2,\ldots,r-2\, .
$$
In particular, $\epsilon_2\le q_0$. On the other hand, by Lemma  
\ref{l3.5}(2) we must have $m_{r-2}\le q+q_0$ and so, by
(\ref{eq3.6}), we find that $\epsilon_2\ge q_0$, i.e. $\epsilon_2=q_0$. 

Finally we show that $r=4$. $\epsilon_2=q_0$ 
implies $\epsilon_{r-2}\le q_0$. Since $m_2\le q+q_0$ (cf. Lemma 
\ref{l3.5}(2)), by (\ref{eq3.6}), we have $\epsilon_{r-2}\ge q_0$.
Therefore $\epsilon_{r-2}=q_0=\epsilon_2$, i.e. $r=4$.
\end{proof}
{\bf Proof of Theorem \ref{B}.} Let $P_1\in X(\fq)$ be as above. By 
(\ref{eq3.7}), Lemma \ref{l3.5}(3) and Lemma \ref{l3.6} we have the
following equation   
$$
y_2^q-y_2=g_2^{q_0}(x^q-x)\ ,
$$
where $g_2$ has no pole except at $P_1$. Moreover, by (\ref{eq3.6}), 
$m_2=q_0+q$ and so $v(g_2)=-q$ (cf. Lemma \ref{l3.5}(3)). Thus 
$g_2=ax+b$ with $a,b\in \fq$, $a\neq 0$, and after some $\fq$-linear
transformations we obtain Theorem \ref{B}.
\begin{remarks*} (1) From the above computations we conclude that the
Deligne-Lusztig curve associated to the Suzuki group $X$ is equipped with
a complete simple base-point-free $g^4_{q+2q_)+1}$, namely
$\cD=|(q+2q_0+1)P_0|$, $P_0\in X(\fq)$. Such a linear system is an
$\fq$-invariant. The orders of $\cD$ (resp. the $\fq$-Frobenius orders)
are $0, 1, q_0, 2q_0$ and $q$ (resp. $0, q_0, 2q_0$ and $q$).

(2) There exists $P_1\in X(\fq)$ such that the $(\cD,P_1)$-orders are
$0,1,q_0+1, 2q_0+1$ and $q+2q_0+1$ (Lemma \ref{l3.3}). Now 
we show that the above sequence is, in fact, the
$(\cD,P)$-orders for each $P\in X(\fq)$. To see this, notice that
$$
\deg(S)=(3q_0+q)(2g-2)+(q+4)(q+2q_0+1)=(4+2q_0)\#X(\fq).
$$
Let $P\in X(\fq)$. By (\ref{eq3.5}), we conclude that 
$v_P(S^\cD)=\sum_{i=1}^{4}(j_i(P)-\nu_{i-1})=4+2q_0$ and so, by
(\ref{eq1.4}), 
that $j_1(P)=1$, $j_2(P)=q_0+1$, $j_3(P)=2q_0+1$, and $j_4(P)=q+2q_0+1$.

(3) Then, by (\ref{eq3.2}) $H(P)$, $P\in X(\fq)$, contains the 
semigroup\newline $H:= \langle q,q+q_0,q+2q_0,q+2q_0+1\rangle$. Indeed,  
$H(P)=H$ since $\#(\mathbb N\setminus H)=g=q_0(q-1)$ (this can be 
proved as in the remark after Lemma \ref{l3.4}; see also 
\cite[Appendix]{hsti}). 

(4) We have  
$$
\deg(R^\cD)=\sum_{i=0}^{4}\epsilon_i(2g-2)+5(q+2q_0+1)=(2q_0+3)\#X(\fq)\,
,
$$
for $P\in X(\fq)$, $v_P(R^\cD)=2q_0+3$ as follows from items (1), (2) and
(\ref{eq1.1}). Therefore the set of $\cD$-Weierstrass points of $X$ is
equal to $X(\fq)$. In particular, the $(\cD,P)$-orders for $P\not\in
X(\fq)$ are $0, 1, q_0, 2q_0$ and $q$. 

(5) We can use the above computations to obtain information on orders
for the canonical morphism. By using the fact that $(2q_0-2)\cD$ is
canonical (cf. Lemma 
\ref{l3.5}(1)) and item (4), we see that the set $
\{a+q_0b+2q_0c+qd: a+b+c+d \le 2q_0-2\}
$
is contained in the set of orders for $\cK_X$ at non-rational points. (By
considering first order differentials on $X$, similar computations were 
obtained in \cite[\S4]{gsti}.)

(6) Finally, we remark that $X$ is 
non-classical for the canonical morphism: we have two different proofs for
this fact: loc. cit. and \cite[Prop. 1.8]{fgt}).
\end{remarks*}
\begin{center}
{\bf Appendix:} A remark on the Suzuki-Tits ovoid
\end{center}
\smallskip

For $s\in \N$, let $q_0:=2^s$ and $q:=2q_0$. It is well known that the
Suzuki-Tits ovoid $\cO$ can be represented in $\P^4(\fq)$ as
$$
\cO=\{(1:a:b:f(a,b):af(a,b)+b^2): a, b \in \fq\}\cup\{(0:0:0:0:0:1)\},
$$
where $f(a,b):=a^{2q_0+1}+b^{2q_0}$ (see \cite{tits}, \cite[p.3]{pent}) 

Let $X$ be the Deligne-Lusztig curve associated to $Sz(q)$ and 
$\cD=|(q+2q_0+1)P_0|$, $P_0\in X(\fq)$ (see \S\ref{3}). By the Remark
(item 3) in \S\ref{3}, 
we can associate to $\cD$ a morphism $\pi=(1:x:y:z:w)$ whose coordinates
satisfy $\dinf(x)=qP_0$, $\dinf(y)=(q+q_0)P_0$, $\dinf(z)=(q+2q_0)P_0$ and
$\dinf(w)=q+2q_0+1$. 
\begin{claim*} (A. Cossidente)\quad 
$\cO=\pi(X(\fq))$.
\end{claim*}
\begin{proof} We have $\pi(P_0)=(0:0:0:0:1)$; we can choose 
$x$ and $y$ satisfying \newline $y^q-y=x^{q_0}(x^q-x)$, 
$z:= x^{2q_0+1}+y^{2q_0}$, and $w:=
xy^{2q_0}+z^{2q_0}=xy^{2q_0}+x^{2q+2q_0}+y^{2q}$
(cf. \cite[\S1.7]{hsti}). For $P\in X(\fq)\setminus\{P_0\}$ set $a:=x(P)$, 
$b:=y(P)$, and $f(a,b):= z(a,b)$. Then $w(a,b)=af(a,b)+b^2$ and we are
done.
\end{proof}
\begin{remark*} The morphism $\pi$ is an embedding. Indeed, since
$j_1(P)=1$ for each $P$ (cf. Remarks \S3(2)(4)), it is enough to see that
$\pi$ is injective. By (\ref{eq3.1}), the points $P$ where $\pi$ could not
be injective satisfy: $P\not\in X(\fq)$$, \frx^3(P)=P$ or $\frx^2(P)=P$.
Now from the Zeta function of $X$ one sees that $\#X(\mathbb
F_{q^3})=\#X(\mathbb F_{q^2})=\#X(\fq)$, and the remark follows.
\end{remark*}
\begin{remark*} From the claim, (\ref{eq3.1}) and \cite{he} we have
$$
{\rm Aut}_{\bar\fq}(X)={\rm Aut}_{\fq}(X)\cong
\{A\in PGL(5,q): A\cO=\cO\}\, .
$$
\end{remark*}

\end{document}